\documentstyle[12pt]{article}
\newcommand{\Sub}{\mathop{\rm Sub}\limits}

\newcommand{\const}{\mathop{\rm const}\limits}

\begin{document}
 \begin{center}

{\bf GRAND LEBESGUE NORM  ESTIMATION \\

\vspace{4mm}

FOR BINARY RANDOM VARIABLES,\\

\vspace{4mm}

 with applications.}\\

\vspace{4mm}

{\sc Eugene Ostrovsky, Leonid Sirota}\\

\vspace{3mm}

 Bar-Ilan University,  59200, Ramat Gan, ISRAEL; \\

 e-mail: eugostrovsky@list.ru \\
 e-mail: sirota3@bezeqint.net \\

 \vspace{4mm}

        {\bf Abstract}

\end{center}

\vspace{3mm}

 \  We calculate the so-called {\it Rademacher's }
 Grand Lebesgue Space norm for a centered  (shifted) indicator (Bernoulli's, binary)  random variable.\par

 \ This norm is optimal for the centered and bounded random variables (r.v.) \par

 \  Using this result we derive a very simple bilateral sharp exponential
 tail estimates for sums of these variables, not  necessary to be
 identical distributed,  under non-standard norming,
  and give some examples to show  the exactness of our estimates. \par

\vspace{3mm}


\vspace{3mm}

 \section { Introduction. Notations. Statement of problem}

   \ In order to formulate our result, we need to introduce some
notations and conditions.  Let $  \{\Omega, B, {\bf P}  \} $ be certain non-trivial probability space.
Let also $ \phi = \phi(\lambda), \lambda \in (-\lambda_0, \lambda_0), \
\lambda_0 = \const \in (0, \infty] $ be some even strong convex which takes positive values for positive arguments twice continuous
differentiable function, such that

$$
 \phi(0) = 0, \ \phi^{//}(0) > 0, \  \exists \lim_{\lambda \to \lambda_0} \phi(\lambda)/\lambda > 0, \eqno(1.1)
$$
including the case when the last limit is equal to infinity. \par

  \ We denote the set of all these function as $ \Phi; \ \Phi = \{ \phi(\cdot) \}. $  \par

 \ We say by definition that the {\it centered} (mean zero) random variable (r.v) $ \xi = \xi(\omega) $
belongs to the Banach space $ B(\phi), $ if there exists some non-negative constant $ \tau \ge 0 $ such that

$$
\forall \lambda \in (-\lambda_0, \lambda_0) \ \Rightarrow
{\bf E} \exp(\lambda \xi) \le \exp[ \phi(\lambda \ \tau) ]. \eqno(1.2).
$$

 \ These spaces appears at first in the article  \cite{Kozatchenko1}.
The complete description and investigation of such a spaces may be found in a monograph \cite{Ostrovsky1},
chapter 1, section 1, p. 22-24. In particular, it was proved that these spaces are really complete Banach spaces. \par

 \ The function $  \lambda \to {\bf E} \exp(\lambda \xi)  $ is said to be  a
moment generating function  for the r.v. $  \xi, $
if there exists at least for one non  zero value $  \lambda. $ \par

 \ The  value $  \lambda_0 $ in the considered in this report examples will be equal to infinity:
$  \ \lambda_0 = \infty. $ \par

 \vspace{3mm}

{\bf Example 1.1.} Let $ \eta $ be a  (renormed) Rademacher's random variable:

$$
{\bf P} (\eta = 1/2) = {\bf P} (\eta = - 1/2)  = 1/2;
$$
then

$$
{\bf E} \exp(\lambda \eta) = \cosh (\lambda/2).
$$

 \ Denote $  \phi_R(\lambda) = \ln \cosh (\lambda/2);  $ then $  \eta \in B \phi_R $ and $ ||\eta| B\phi_R = 1. $ \par

 \ The function  $  \phi_R(\lambda) = \ln \cosh (\lambda/2)  $  one can named as a {\it natural function} for the Rademacher's
random variable; see exact definition further. \par

\vspace{3mm}

 \ Evidently,

$$
 \phi_R(\lambda)   \sim  \lambda^2/8, \ \lambda \to 0; \   \phi_R(\lambda)   \sim  |\lambda|/2, \ |\lambda|  \to \infty.
$$

 \vspace{3mm}

 \ The minimal value $ \tau $ satisfying (1.2) is called a $ B(\phi) \ $ norm
of the variable $ \xi, $ write

 $$
 ||\xi||B(\phi) = \inf \{ \tau, \ \tau > 0: \ \forall \lambda \ \Rightarrow
 {\bf E}\exp(\lambda \xi) \le \exp(\phi(\lambda \ \tau)) \}. \eqno(1.3)
 $$

 \  The correspondent Grand Lebesgue Space norm  $ ||\xi||G\phi_R $ will be named
{\it Rademacher's } norm. \par

\vspace{4mm}

 \  These spaces are very convenient for the investigation of the r.v. having a
exponential decreasing tail of distribution, for instance, for investigation of the limit theorem,
the exponential bounds of distribution for sums of random variables,
other non-asymptotical properties of the random vectors and processes, problem of continuous of random fields,
study of Central Limit Theorem in the Banach space etc.\par

 \  The space $ B(\phi) $ with respect to the norm $ || \cdot ||B(\phi) $ and
ordinary algebraic operations is a Banach space which is isomorphic to the subspace
consisted on all the centered variables of Orlicz's space $ (\Omega,F,{\bf P}), N(\cdot) $ with $ N \ - $ function

$$
N(u) = \exp(\phi^*(u)) - 1, \ \phi^*(u) = \sup_{\lambda} (\lambda u -
\phi(\lambda)). \eqno(1.4)
$$

 \ The transform $ \phi \to \phi^* $ is called Young-Fenchel, or Legendre transform. The proof of considered
assertion used the properties of saddle-point method and theorem of Fenchel-Moraux:
$$
\phi^{**} = \phi.
$$

 \ The next facts about the $ B(\phi) $ spaces are proved in \cite{Kozatchenko1}, \cite{Ostrovsky1}, p. 19 - 40.

$$
{\bf 1.} \ \xi \in B(\phi) \Leftrightarrow {\bf E } \xi = 0, \ {\bf and} \ \exists C = \const > 0,
$$

$$
U(\xi,x) \le \exp(-\phi^*(Cx)), x \ge 0, \eqno(1.5)
$$
where $ U(\xi,x)$ denotes in this article the {\it tail} of
distribution of the r.v. $ \xi: $

$$
U(\xi,x) = \max \left( {\bf P}(\xi > x), \ {\bf P}(\xi < - x) \right),
\ x \ge 0,
$$
and this estimation is in general case asymptotically exact. \par

 \ Here and further $ C, C_j, C(i) $ will denote the non-essentially positive
finite "constructive" constants.\par

 \ More exactly, if $ \lambda_0 = \infty, $ then the following implication holds:

$$
\lim_{\lambda \to \infty} \phi^{-1}(\log {\bf E} \exp(\lambda \xi))/\lambda =
K \in (0, \infty) \eqno(1.6a)
$$
if and only if

$$
\lim_{x \to \infty} (\phi^*)^{-1}( |\log U(\xi,x)| )/x = 1/K, \eqno(1.6b)
$$
 see \cite{Bagdasarov1}. \par

 \ Here and further $ f^{-1}(\cdot) $ denotes the inverse function to the
function $ f $ on the left-side half-line $ (C, \infty). $ \par

 \ The function $ \phi(\cdot) $ may be "constructive" introduced by the formula
$$
\phi(\lambda) = \phi_0(\lambda) \stackrel{def}{=} \log \sup_{t \in T}
 {\bf E} \exp(\lambda \xi(t)), \eqno(1.7)
$$
 if obviously the family of the centered r.v. $ \{ \xi(t), \ t \in T \} $ satisfies the {\it uniform } Kramer's condition:
$$
\exists \mu \in (0, \infty), \ \sup_{t \in T} U(\xi(t), \ x) \le \exp(-\mu \ x),
\ x \ge 0. \eqno(1.8)
$$

 \ In this case, i.e. in the case the choice the function $ \phi(\cdot) $ by the
formula (1.7), we will call the function $ \phi(\lambda) = \phi_0(\lambda) $
a {\it natural } function for the family of the centered r.v. $ \{ \xi(t), \ t \in T \}. $ \par

\vspace{3mm}

  \ We say that the {\it centered:} $ {\bf E} \xi = 0 $ numerical random variable (r.v.)
 $ \xi = \xi(\omega), \ \omega \in \Omega $ is subgaussian, or equally, belongs to the space $ \Sub(\Omega), $
 if there exists some non-negative constant $ \tau \ge 0 $ such that

$$
\forall \lambda \in R  \ \Rightarrow
{\bf E} \exp(\lambda \xi) \le \exp[ \lambda^2 \ \tau^2 ]. \eqno(1.9)
$$
 The minimal value $ \tau $ satisfying (1.1) is called a  subgaussian  norm
of the variable $ \xi, $ write

 $$
 ||\xi||\Sub = \inf \{ \tau, \ \tau > 0: \ \forall \lambda \in R \ \Rightarrow {\bf E}\exp(\lambda \xi) \le
 \exp(\lambda^2 \ \tau^2) \}.
 $$

 Evidently,

$$
||\xi||\Sub = \sup_{\lambda \ne 0} \left[ \sqrt{ \ln {\bf E}  \exp ( \lambda \xi)  }/|\lambda| \right].  \eqno(1.10)
$$

 \ So, the space $ \Sub = \Sub(\Omega) $ is the particular case of the general $  B(\phi) $ spaces with
 $ \phi(\lambda) =\phi_2(\lambda) = \lambda^2, \ \lambda \in R.  $ \par

 \ This important notion was introduced before the appearing of the general theory of $  B(\phi) \ $ spaces
  by  J.P.Kahane \cite{Kahane1}; V.V.Buldygin and Yu.V.Kozatchenko \cite{Buldygin1} proved
that the set $ \Sub(\Omega) $  relative the norm $  ||\cdot|| $ is complete Banach space which is isomorphic to subspace
consisting only from the centered variables of Orlicz's space over $ (\Omega, B,P)  $ with $ N \ - $ Orlicz-Young function
 $ N(u) = \exp(u^2) - 1 $  \cite{Kozatchenko1}.  \par

   If $ ||\xi||\Sub = \tau \in (0,\infty),  $ then

 $$
 \max [{\bf P}(\xi > x),  {\bf P}(\xi < -x)  ] \le \exp(- x^2/(4 \tau^2)  ), \ x \ge 0;
 $$
  and  the last inequality is in general case non-improvable.  It is sufficient for this to consider the case when
 the r.v. $  \xi  $ has the centered Gaussian non-degenerate distribution.\par

  Conversely, if  $ {\bf E} \xi = 0 $ and if  for some positive finite constant $  K  $

 $$
 \max [{\bf P}(\xi > x),  {\bf P}(\xi < -x)  ] \le \exp(- x^2/K^2  ), \ x \ge 0,
 $$
 then $ \xi \in \Sub(\Omega) $ and $ ||\xi||\Sub < 4 K. $ \par

 The subgaussian norm in the subspace of the centered r.v. is equivalent to the following Grand Lebesgue Space (GLS)
 norm:

 $$
|||\xi||| := \sup_{s \ge 1} \left[ \frac{|\xi|_s}{\sqrt{s}} \right], \hspace{6mm}
|\xi|_s \stackrel{def}{=}  \left[ {\bf E} |\xi|^s \right]^{1/s}.
 $$

 \ For the non - centered r.v. $ \xi $  the subgaussian norm may be defined as follows:

 $$
 ||\xi|| \Sub := \left[  \left\{ ||\xi - {\bf E} \xi||\Sub \right\}^2 + ( {\bf E} \xi)^2  \right]^{1/2}.
  $$

  \ More detail investigation of these spaces see in the monograph \cite{Ostrovsky1}, chapter 1.  \par

 \ We denote as usually by $  I(A) = I(A; \omega), \ \omega \in \Omega, \ A \in B  $ the indicator function of event $  A. $
 Further, let $  p  $ be arbitrary number from the set $ [0,1]: \ 0 < p < 1 $ and let $ A(p) $  be any event such that
 $ {\bf P}(A(p)) = p. $ Denote also  $ \eta_p = I(A(p)) - p; $  the {\it  centering  } of the r.v. $ I(A(p)); $
 then $ {\bf E} \eta_p = 0  $ and

 $$
 {\bf P}(\eta_p = 1 - p) = p; \hspace{6mm}  {\bf P}(\eta_p = -p) = 1 - p.   \eqno(1.11)
 $$

 \ The case $ p = 1 - p = 1/2  $ correspondent to the considered before  case of Rademacher's random variable. \par

\vspace{3mm}

{\bf  Our goal in this short report is to  investigate the value of the Rademacher's norm for the random variable $ \eta_p. $

\vspace{3mm}

 \  We derive in the third section a very simple  non-asymptotical bilateral
 tail estimates for sums of these variables, not  necessary to be
 identical distributed,  under non-standard norming.} \par

\vspace{3mm}

 \  Let us describe briefly some previous works.
 Define the following non-negative continuous on the closed segment $ p \in [0,1]  $ function

$$
Q(p) =  \sqrt{ \frac{1 - 2p}{4 \ln(( 1 - p )/p)} },\eqno(1.12)
$$
so that $ Q(0+0) = Q(1-0) = 0 $  and $ Q^2(1/2) = 1/8 $ (Hospital's rule).  Note also

$$
p \to 0+ \ \Rightarrow Q(p) \sim \frac{0.5}{\sqrt{|\ln p|}}, \hspace{6mm}  p \to 1 - 0 \Rightarrow  Q(p) \sim \frac{0.5}{\sqrt{|\ln(1 - p)|}}.
\eqno(1.13)
$$

 \ The last circumstance play a very important role in the non-parametrical statistics,  see \cite{Gaivoronsky1}, \cite{Kiefer1}.\par

\vspace{3mm}

 It is known  \cite{Kearns1}, \cite{Berend1}, \cite{Schlemm1}, \cite{Buldygin2}, \cite{Ostrovsky101}, \cite{Ostrovsky102}
that

 $$
 ||\eta_p||\Sub = Q(p).
 $$

\vspace{3mm}

 \ Applications of these estimates in the non-parametrical statistics may be found in the articles
\cite{Gaivoronsky1}, \cite{Kiefer1}. Other  application is described in \cite{Chen1}. \par

 \ Another approach  and applications see in the works
\cite{Bentkus1}, \cite{Berend1},  \cite{Pinelis1}, \cite{Raginsky1}, \cite{Schlemm1}, \cite{Serov1}, \cite{Zubkov1} etc. \par

\vspace{3mm}

 \section{ Auxiliary result.}

 \vspace{3mm}

 \hspace{3mm}  Recall first of all that

$$
y = y(z):= \cosh^{-1} z = \ln (z \pm \sqrt{z^2 - 1}), \ z \ge1.
$$

 \ We agree to take only the following branch of these function

$$
y(z) = \cosh^{-1} z = \ln (z + \sqrt{z^2 - 1}), \ z \ge 1.
$$

 \ Note that

$$
z \to 1 + 0 \ \Rightarrow y(z) \sim \sqrt{ 2(z-1) }, \eqno(2.0a)
$$

$$
z \to \infty \ \Rightarrow y(z) \sim \ln z. \eqno(2.0b)
$$

 \ The natural function for the family of the (centered) r.v. $  \{  \eta_r \},  \ 0 < r < 1  $ has a form

$$
 \beta_r(\lambda) \stackrel{def}{=} {\bf E} e^{\lambda \eta_r} = r e^{\lambda (1 - r) } + (1 - r) e^{ - r \lambda }, \
 \lambda \in (-\infty, \infty), \ 0 < r < 1, \eqno(2.1)
$$
so that

$$
\beta_{1/2}(\lambda) = 0.5 \left( e^{\lambda/2} + e^{-\lambda/2} \right) =\cosh(\lambda/2).
$$

 \ Evidently,

$$
\lambda \to \infty \ \Rightarrow  \beta_r(\lambda) \sim r e^{\lambda(1-r)}, \  r = \const \in (1/2,1),
$$

$$
\lambda \to 0 \ \Rightarrow  \beta_r(\lambda) \sim 1 + 0.5 \lambda^2 r (1 - r), \  r = \const \in (0,1),
$$

\vspace{3mm}

 \ Introduce an important function, which may be named as Rademacher's norm of the binary random variable,

$$
g_R(r) = g(r) \stackrel{def}{=}  \sup_{\lambda \ne 0} \left[ \frac{\cosh^{-1} [\beta_r(\lambda)] } {|\lambda|/2}  \right]=
$$

$$
 \sup_{\lambda \ne 0} \left[ \frac{\cosh^{-1}( r e^{\lambda (1 - r) } + (1 - r) e^{ - r \lambda }  )}{|\lambda|/2}  \right],
 \ r \in (0,1). \eqno(2.2)
$$

\vspace{3mm}

{\bf Proposition 2.1.} It follows immediately from the direct definition of the $ B(\phi_R)  $ norm that

$$
||\eta_r||B\phi_R = g(r), \hspace{5mm} 0 < r < 1. \eqno(2.3)
$$

\vspace{3mm}

 \ Let us itemize now some important for us properties of introduced function $ g = g(r) = g_R(r), \ 0 \le r \le 1. $
All this properties may be easily obtained from the known asymptotical behavior of both the functions $  y(z) $
and $  \beta_r(\lambda).  $  \par

 \vspace{3mm}

 {\bf 1.}   This function is bounded and continuous on the closed interval $ [0,1]. $
 More detail: the inequality $ 0 < g(r)  \le 2  $ is obvious. \par

  \ Moreover

 $$
 g(0+) = g(1-0) = 2.
 $$

 \ Note that the last equality stand in contradiction to the analogous fact (1.13) for the subgaussian norm for at
the same binary random  variable. \par

 \ As a consequence: the function $  g = g(r)  $ can be defined as a continuous positive function on the closed
 interval $  [0,1] $   such that $  g(0) = g(1) = 2.  $\par

 \  So, $  \max_{r \in [0,1]} g(r) = 2 = g(0) = g(1).  $ \par

\vspace{3mm}

{\bf 2.}  On the other hand,  we obtain  after some calculations

$$
g(r) \ge \lim_{\lambda \to 0} \left[ \frac{\cosh^{-1} [\beta_r(\lambda)] } {|\lambda|/2}  \right] = 2 \sqrt{r(1 - r)}, \
r \in (0,1) \eqno(2.4)
$$

 \vspace{3mm}

{\bf 3.} \ Evidently, $  g(1-r) = g(r), $ (symmetry), so that it is enough to investigate this function only on the interval
$  1/2 \le r \le 1.  $ \par

 \vspace{3mm}

{\bf 4.}  It is easy to calculate $  g(1/2) = 1. $ \par

 \vspace{3mm}

{\bf 5.}    Note in addition

$$
g(r) \ge \lim_{|\lambda| \to \infty} \left[ \frac{\cosh^{-1} [\beta_r(\lambda)] } {|\lambda|/2}  \right] =
2 \max(r, 1-r), \ r \in (0,1),
$$
but the last function is less than  $ 2 \sqrt{r(1 - r)}. $ \par

\vspace{4mm}

 \ The following rough estimate will be practically used in the next section. \par

\vspace{3mm}

{\bf Proposition 2.2.} \\

\vspace{2mm}

$$
  \sup_{r \in [0,1]} \ ||\eta_r||G\psi_R = 1, \eqno(2.5)
$$
 or equally

 $$
\sup_{r \in [0,1]} {\bf E} e^{\lambda \eta_r} = \cosh (\lambda/2). \eqno(2.5a)
 $$

\vspace{3mm}

{\bf Proof.} \\

\vspace{2mm}

{\bf 1.} It is sufficient to consider for reasons of symmetry only the cases $  r \in [1/2, 1]  $ and
analogously $ \lambda \ge 0. $\\

\vspace{3mm}

{\bf 2.}  Further, we have proved the following equivalent elementary inequality

$$
 \beta_r(\lambda) =
r e^{\lambda (1 - r)} + (1 - r) e^{-\lambda r} \le \cosh(\lambda/2), \ r \in (1/2, 1), \eqno(2.6)
$$
wherein $  \lambda \ge 0; $ the cases $ r = 1/2, \ r = 1 $ and $  r = 0 $ are trivial.\par

\vspace{3mm}

{\bf 3.} Put  for simplicity  $  r = 1/2 + \delta, \ \delta \in (0,1/2), $ then we deduce after some calculations

$$
\beta = e^{ - \lambda \delta } \left[ \cosh (\lambda/2) + 2 \delta \sinh(\lambda/2)   \right].
$$

\vspace{3mm}

{\bf 4.} Our inequalities (2.5) and (2.6) takes the form

$$
\left[ \cosh (\lambda/2) + 2 \delta \sinh(\lambda/2)   \right] \le \cosh(\lambda/2)
$$
or equally

$$
 2 \delta \sinh(\lambda/2) \le \cosh(\lambda/2) \left[ e^{\lambda \delta} - 1 \right]. \eqno(2.7)
$$

 \ The inequality (2.7) follows in turn taking into account the positivity of the value of product $  \lambda \delta $
from the one of the form

$$
\sinh \mu \le 2 \mu \ \cosh \mu, \hspace{5mm}  \mu = 2 \lambda > 0. \eqno(2.8)
$$

\vspace{3mm}

{\bf 5.} The last inequality (2.8) may be elementary proved by means of juxtapositions of correspondent Taylor's members.\par

\vspace{3mm}

{\bf 6.} The equality in the assertion of proposition  2.2 is reached for example for the value $  r = 1/2 $
as well as as $  \lambda \to \infty $ and as $ r \to 1 - 0 $ or equally as $  r \to 0+. $\par

\vspace{3mm}

 \ This completes the proof of proposition 2.2.\par

\vspace{3mm}

 \section{ Main result: tail estimations for sums of independent indicators under non - standard norming.}

 \vspace{3mm}

  \ Let $  p(i), \ i = 1,2, \ldots,n  $ be positive numbers such that $ 0 < p(i) < 1,  $ and let $ A(i) $ be
 {\it independent} events for which $  {\bf P}(A(i)) = p(i). $  Introduce a sequence of two - values (binary, generalized
Rademacher's independent random variables $  \zeta := \{ \zeta(i) \}, \hspace{5mm}
\zeta(i) = I(A(i)) - p(i), $ and define its sum

$$
 S(n) := w(n)^{-1} \sum_{i=1}^n \zeta(i),\eqno(3.1)
$$
where the norming function $  w = w(n)  $ is any deterministic strictly increasing to infinity numerical sequence
such that  $  w(1) = 1. $ \par

\vspace{3mm}

 \ {\sc We intend in this section to derive the bilateral uniform exponential  bounds for the tail of distribution }

$$
T_w(u) \stackrel{def}{=} \sup_n \sup_{ \{ \zeta\} } \max \left\{ {\bf P} (S(n) > u),  \  {\bf P} (S(n) < - u) \right\},
\ u > 1, \eqno(3.2)
$$
where $ \sup_{ \{ \zeta\} } $ is calculated  over all the centered Rademacher's independent random variables
$  \zeta := \{ \zeta(i) \}. $ \par

 \  The normalization $   w_{1/2}(n) = \sqrt{n}  $ can be considered a classic, see e.g. \cite{Kozatchenko1},
\cite{Ostrovsky1}, chapter 2,  sections 2.1 and 2.2.  For instance, if the r.v. are i., i.d, centered and

$$
\max \left\{ {\bf P} (\zeta(i) > u),  \  {\bf P} (\zeta(i) < - u) \right\} \le \exp \left(-u^k \right), \ u \ge 0,
\ k = \const > 0,
$$
then

$$
T_{w_{1/2}}(u) \le \exp \left( - C(k) \ u^{\min(k,2)} \right), \ 0 < C(k) = \const < \infty,
$$
and the last estimate  is essentially non-improvable. \par
  \ It is known for instance, see \cite{Ostrovsky1}, chapter 1, section 1.6 that

$$
||\sum_{i=1}^n \zeta(i)||\Sub \le \sqrt{ \sum_{i=1}^n (||\zeta(i)||\Sub)^2  }.
$$

\vspace{3mm}

 \ Therefore, it is reasonable to suppose

$$
\lim_{n \to \infty} w(n)/\sqrt{n} = \infty;
$$
other case is trivial for us. \par

 \ On the on the other hand, if

$$
\lim_{n \to \infty} w(n)/n = \infty,
$$
then evidently

$$
T_{w_{1/2}}(u) = 0,  \hspace{5mm} u > 1.
$$

  \ Thus, we must exclude also both these cases. \par

 \vspace{3mm}

  \ The exact formulating using for us assumptions will be specified below.\par

\vspace{3mm}

 \ We shall touch briefly earlier work in the considered here problem. The particular case of our statement of problem, even
for the sequence of martingale differences, may be found in the articles  \cite{Fan1},
\cite{Lesign1}, \cite{Li1}, \cite{Liu1}, \cite{Ostrovsky103}, \cite{Ostrovsky104}. \par

 \ A very interest application to the investigation  of the free energy of directed polymers in random environment
is described in the article belonging to Liu Q. and Watbled F. \cite{Liu1}. \par

\vspace{3mm}

 \ Let us now itemize some conditions imposed on the norming function $  w = w(n).  $ \par

\vspace{3mm}

{\bf A1.} There exists a strictly decreasing twice continuous differentiable function,  defied on the set $ \lambda \ge 1,  $
which we will denote also $  w = w(\lambda), $ such that $  w(\lambda)/_{\lambda = n} = w(n), \ n = 1,2,3, \ldots. $\par

\vspace{3mm}

{\bf A2.}

$$
\lambda \to \infty \ \Rightarrow \frac{w(\lambda)}{\lambda} \downarrow \ 0. \eqno(3.3)
$$

\vspace{3mm}

{\bf A3.}

$$
\lambda \to \infty \ \Rightarrow \frac{w(\lambda)}{ \sqrt{\lambda}} \uparrow \ \infty. \eqno(3.4)
$$

\vspace{3mm}

{\bf A4.} The inverse  function $  \lambda \to w^{-1}(\lambda), \ \lambda \ge 1 $ is convex. \par

\vspace{3mm}

{\bf A5.} The function $  \lambda \to w(\lambda), \ \lambda \ge 1 $ satisfies the $  \Delta_2 \ - $ condition

$$
\sup_{\lambda > 1} \left[\frac{w(2 \lambda)}{w(\lambda)} \right] < \infty. \eqno(3.5)
$$

\vspace{3mm}

 \ Define also  a new function

 $$
 v(u) = v_w(u) :=  \left(w^{-1} \right)^*(u), \ u \ge 1.
 $$

 \ Recall that the transformation $  f \to f^* $ is named Young-Fenchel, or Legendre transform, see (1.4). \par

\vspace{3mm}

 {\bf Theorem 3.1.}   {\it Let all the conditions  {\bf A1} - {\bf A5} be satisfied. We assert then as}
$  u \to \infty $

$$
| \ \ln  T_w(u) \ | \asymp  v_w(u). \eqno(3.6)
$$

\vspace{3mm}

{\bf Proof.} \\

\vspace{3mm}

{\bf 1.} Let  us calculate (and evaluate) first of all the moment generating function  for the sequence of r.v. $  S(n). $
We  have using the independence of the r.v. $  \zeta(i) $

$$
{\bf E}e^{\lambda S(n)} = \prod_{i=1}^n {\bf E} e^{\lambda \zeta(i)/w(n)}  =
\prod_{i=1}^n \beta_{p(i)} \left( \frac{\lambda}{w(n)} \right).
$$

{\bf 2.} One can use the proposition 2.2, more exactly, the estimate (2.5a):

$$
{\bf E}e^{\lambda S(n)} \le \cosh^n \left(  \frac{\lambda}{w(n)}  \right) =
\exp \left[ n \ln \cosh \left( \frac{\lambda}{w(n)} \right)   \right], \eqno(3.7)
$$
wherein the last inequality (3.7) is sharp: is achievable for instance when  all the r.v. $  \zeta(i)  $ are Rademacher's. \par
 On the other words, we can and will suppose all the independent
 variables $ \zeta(i)  $ have the ordinary symmetrized  Rademacher' distribution,\par

\vspace{3mm}

{\bf 3.} Therefore

$$
 \sup_n {\bf E}e^{\lambda S(n)} \le \sup_n \exp \left[ n \ln \cosh \left( \frac{\lambda}{w(n)} \right)   \right], \eqno(3.8)
$$
and it is easily to derive using the known properties of the function $  w(\cdot) $

$$
\ln \sup_n {\bf E}e^{\lambda S(n)} \asymp  w^{-1}(\lambda), \ |\lambda| > 1;  \eqno(3.9)
$$
the case $  |\lambda| \le 1 $ is simple. \par

 \ In particular, if we choose in the  right hand (3.8) $  w(n) = \lambda,  $ then

$$
\ln \sup_n {\bf E}e^{\lambda S(n)} \ge  C_1(w)  \cdot w^{-1}(\lambda), \ |\lambda| > 1.  \eqno(3.9a)
$$

 \ In detail,  let $  \lambda > 1. $ There exists an unique value $  n_0 = n_0(w, \lambda) $ such that

$$
 w^{-1}(\lambda) \le n_0 < w^{-1}(\lambda) + 1.
$$
 Then

$$
 \sup_n {\bf E}e^{\lambda S(n)} \ge  {\bf E}e^{\lambda S(n_0)},
$$
and we have consequently $  n_0 \ge w^{-1}(\lambda);  $

$$
\frac{\lambda}{w(n_0)} \ge \frac{\lambda}{w(w^{-1}(\lambda) +1) } \ge
\frac{\lambda}{ \lambda + w(1) } \ge \frac{1}{ 1 + w(1) },
$$
we exploited the convexity of the function $  w^{-1}(\cdot),  $ condition {\bf A4}. \par

 Thus, one can  to choose in (3.10)

$$
C_1(w) = \frac{1}{ 1 + w(1) }.
$$

\vspace{3mm}

 \ We turn now to the withdrawal of the upper bound for the value $ Z := \sup_n {\bf E}e^{\lambda S(n)}. $ Define an
absolute constant

$$
C = e + 1/e - 2 \approx 1.0862 \ldots.
$$

 For our purpose we estimate:

$$
\ln \cosh \lambda \le \lambda, \ \lambda \ge 1;
$$

$$
 |\lambda| < 1 \ \Rightarrow \hspace{4mm}
\cosh \lambda =  1 + \frac{\lambda^2}{2!} +  \frac{\lambda^4}{4!} + \frac{\lambda^6}{6!} + \ldots \le
$$

$$
1 +  \frac{\lambda^2}{2} \times \left[ 1 + 2 \cdot \left( \frac{1}{4!} +  \frac{1}{6!} + \ldots \right) \right] =
$$

$$
1 +  \frac{\lambda^2}{2} \times \left[ 1 + 2 \cdot (\cosh 1 - 3/2) \right] = 1 +  C \cdot \frac{\lambda^2}{2};
$$

 $$
|\lambda| < 1 \ \Rightarrow \hspace{4mm} \ln \cosh \lambda \le C \cdot \frac{\lambda^2}{2}.
 $$

 \ We find combining the obtained estimates for the positive values $  \lambda: $

$$
\ln Z = \ln \sup_n {\bf E}e^{\lambda S(n)} \le \lambda \cdot \frac{n}{w(n)} \cdot I(n \le w^{-1}(\lambda)) +
$$

$$
C \cdot \frac{\lambda^2 n}{w^2(\lambda)} \cdot I(n \ge w^{-1}(\lambda)),
$$
where $  I(A) $ denotes the indicator function for the predicate $  A. $ \par

 \ Obviously,

$$
\ln Z \le \lambda \cdot \frac{\mu}{w(\mu)} \cdot I(1 \le \mu < w^{-1}(\lambda)) +
$$

$$
C \cdot \frac{\lambda^2 \mu}{w^2(\lambda)} \cdot I(\mu \ge w^{-1}(\lambda)).
$$

 \ It follows immediately from the assumptions {\bf A1} - {\bf A5} that the function of the variable $  \mu, \mu \ge 1 $
in the right-hand side of the last inequality achieved its maximal value at the point $  \mu = w^{-1}(\lambda)  $ and
herewith

$$
Z \le \exp \left(C  w^{-1}(\lambda) \right), \ |\lambda| > 1.
$$

 \ Totally, we obtained the following uniform bilateral estimates for the moment generating function of the random sequence  $  S(n): $

$$
\exp \left(C_1(w) \  w^{-1}(\lambda) \right) \le  \sup_n {\bf E} e^{\lambda S(n)}
\le \exp \left(C  w^{-1}(\lambda) \right), \ |\lambda| > 1. \eqno(3.10)
$$

\vspace{3mm}

{\bf 4.} The proposition of theorem 3.1 follows now from (3.10) and from the main result of the article
\cite{Bagdasarov1}, see also \cite{Ostrovsky1}, chapter 1, section 1.4.\par

\vspace{3mm}

{\bf  Example 3.1.} Let for instance $  w(\lambda) = \lambda^{3/4}, \ \lambda \ge 1;  $ then

$$
\ln \left| \ \sup_n {\bf P} \left( n^{-3/4} \sum_{i=1}^n \zeta(i) > u   \right) \ \right| \asymp  u^4, \ u \ge 1. \eqno(3.11)
$$

 \vspace{3mm}

  \ {\bf Remark 3.1.} Let us introduce the following function

$$
 \theta(\lambda) := \sup_n  \left[ \ln \cosh^n \left(  \frac{\lambda}{w(n)}  \right) \right];
\eqno(3.12)
$$
so that

$$
\sup_n {\bf E}e^{\lambda S(n)} \le e^{ \theta(\lambda)}. \eqno(3.13)
$$

 \ It follows from the inequality (3.13) by means of Chernov's inequality only unilateral inequality

$$
\sup_n \max ( {\bf P}(S(n) > u), {\bf P}(S(n) < - u) ) \le e^{- \theta^*(u)}, \ u \ge 0, \eqno(3.14)
$$
still without all the conditions {\bf A1} -  {\bf A5}.\par

\vspace{3mm}

\section{ Concluding remarks.}

\vspace{3mm}

 \ {\bf A.} It is known, see \cite{Rivasplata1}, that (after commensuration) if $  X $ be a mean zero r.v. $ \ X: \ {\bf E X = 0} \ $
and is bounded a.e.:  $  |X| \le 1/2, $ then

$$
{\bf E} e^{\lambda X} \le \cosh(\lambda/2), \ \lambda \in R.
$$

 \ On the other  words,

$$
|| \ X \ ||B\phi_R \le 1,
$$
herewith the equality in the last estimate is achieved only in the case when $  X  $ has the (symmetrical) Rademacher's
binary distribution. \par

 \ Therefore, all the results of theorem 3.1 remains true for the arbitrary sequence of independent centered such a variables,
not necessary be identical distributed.  \par
 \ This proposition  can be considered as some  complement to the classical theorem  of W.Hoeffding  \cite{Hoeffding1},
see also \cite{Bentkus1}. \par

\vspace{3mm}

 \ {\bf B.} The case of sums of {\it weakly dependent} binary r.v., including sums of martingale differences,
is investigated in the recent article \cite{Pelekis1}.\par

\end{document}